\documentclass[12pt]{article}   	
\usepackage{geometry}                		
\usepackage[parfill]{parskip}    		
\usepackage{graphicx}				
\usepackage[pdftex,pagebackref,letterpaper=true,colorlinks=true,pdfpagemode=none,urlcolor=blue,linkcolor=blue,citecolor=blue,pdfstartview=FitH]{hyperref}
\usepackage{amsmath,amsfonts}
\usepackage{color}
\usepackage{amsthm}
\usepackage{boolexpr}
\usepackage{amssymb,latexsym,pstricks,pst-plot}
\usepackage[english]{babel}
\usepackage{pgf}
\usepackage{tikz}
\usetikzlibrary{arrows,automata}
\usepackage[latin1]{inputenc}
\usepackage[UKenglish]{isodate}
\allowdisplaybreaks

\usepackage{amssymb}
\newtheorem{thm}{Theorem}

\newtheorem{observation}[thm]{Observation}

\title{Structure and asymptotics for Motzkin numbers modulo small primes using automata}
\author{Rob Burns}

\begin{document}
\maketitle
\begin{abstract}
We establish the asymptotic density and some structure results for the Motzkin numbers modulo primes up to $29$ using automata.
\end{abstract}

\section{Introduction}
The Motzkin numbers $M_n$ are defined by
$$
M_n := \sum_{k \geq 0} \binom {n}{2k} C_{k}
$$
where $C_{k}\, $ are the Catalan numbers. 

There has been some work in recent years on analysing the Motzkin numbers $M_n$ modulo primes and prime powers. Deutsch and Sagan \cite{Sagan2006} provided a characterisation of Motzkin numbers divisible by $2, 4$ and $5$. They also provided a complete characterisation of the Motzkin numbers modulo $3$ and showed that no Motzkin number is divisible by $8$. Eu, Liu and Yeh \cite{Eu2008} reproved some of these results and extended them to include criteria for when $M_n$ is congruent to $\{ 2, 4, 6 \} \mod 8$. Krattenthaler and M{\"u}ller \cite{KM2013} established identities for the Motzkin numbers modulo higher powers of 3 which include the modulo 3 result of  \cite{Sagan2006} as a special case. Krattenthaler and M{\"u}ller \cite{Krat2016} have more recently extended this work to a full characterisation of \mbox{$M_n \mod 8$} in terms of the binary expansion of $n$. Their characterisation is rather elaborate and less susceptible to analysis than that provided in \cite{Eu2008}. The results in \cite{KM2013} and \cite{Krat2016} are obtained by expressing the generating function of $M_n$ as a polynomial involving a special function. Rowland and Yassawi \cite{RY2013} investigated $M_n$ in the general setting of automatic sequences.  The values of $M_n$ (as well as other sequences) modulo prime powers can be computed via automata. Rowland and Yassawi provided algorithms for creating the relevant automata. They established results for $M_n$ modulo small prime powers, including a full characterisation of $M_n$ modulo 8 (modulo $5^2$ and $13^2$ are available from Rowland's website). They also established that $0$ is a forbidden residue for $M_n$ modulo $8$, $5^2$ and $13^2$. In theory the automata can be constructed for any prime power but computing power and memory quickly becomes a barrier. For example, the automata for $M_n$ modulo $13^2$ has over 2000 states. Rowland and Yassawi also went on to describe a method for obtaining asymptotic densities of $M_n$.

We will use Rowland and Yassawi's automata to establish asymptotic densities of $M_n$ modulo primes up to $29$. We will also make note of some structure results that appear from an examination of the relevant state diagrams of the automata. We have stopped at $29$ but the method can be applied to any prime power taking into account memory and processing restrictions. We do not provide full details of the automata that we have constructed. Full details are available from our website. The reason for this is that the number of transitions that would need to be included (1015 transitions when looking at $M_n \mod 29$ for example) would make any state diagram unreadable.

Table~\ref{results} summarises the results that will be presented in subsequent sections.

Here, the asymptotic density of a subset $S$ of $\mathbb{N}$ is defined to be
$$
\lim_{N \to \infty} \frac {1}{N} \#\{ n \in S : n \leq N \}
$$
if the limit exists, where $\#S$ is the number of elements in a set $S$. 

For a number $p$, we write the base $p$ expansion of a number $n$ as 
$$
[\, n ]\,_{p} = \langle n_{r} n_{r-1} ... n_{1} n_{0} \rangle
$$ 
where $n_{i} \in [\, 0, p - 1]\,$  and 
$$
n = n_{r} p^{r} + n_{r-1} p^{r-1} + ... + n_{1} p + n_{0}.
$$ 
 
\bigskip

\begin{table}[tbp]
  \centering
   \begin{tabular}{ | c | c | p{10cm} |}
 \hline
 Prime & Density & Other relations \\ \hline
 & & \\
7 & 1 & \mbox{$M_{n} \equiv 0 \mod 7 \,$ when $n = 7^k - 2$.} \\
& & \mbox{$M_{n} \equiv 2 \mod 7 \,$ when $n = 7^k - 1$.} \\ \hline
& & \\
11 & $\frac{1}{55} = \frac{2}{10*11}$ & \\
& & \\ \hline
& & \\
13 & $\frac{1}{78} = \frac{2}{12*13}$ & \\ 
& & \\ \hline
& & \\
& & \mbox{$M_n \equiv 0 \mod 17 \, $ when $\, n = 17^{2k} - 2 \,$.} \\
17 & 1 & \mbox{$M_n \equiv -1 \mod 17 \,$ when $\, n = 17^{2k - 1} - 2 \,$.} \\
& & \mbox{$M_n \equiv 2 \mod 17 \,$ when  $n = 17^{2k} - 1$.} \\
& & \mbox{$M_n \equiv -1 \mod 17 \,$ when $\, n = 17^{2k-1} - 1 \,$.} \\ \hline
& & \\
19 & 1 & \mbox{$M_n \equiv 0 \mod 19 \,$ when $\, n = 19^k - 2 \,$.} \\
& & \mbox{$M_n \equiv 2 \mod 19 \,$ when $ \, n = 19^k - 1 \,$.} \\ \hline
& & \\
23 & $\frac{1}{253} = \frac{2}{22*23}$ & \\
& & \\ \hline
& & \\
29 & $\frac{22}{3045} > \frac{2}{28*29}$ & \\
& & \\ \hline
  \end{tabular}
  \caption{Table of results}
  \label{results}
\end{table}

\bigskip

\section{Background on Motzkin numbers modulo primes}

As mentioned in the introduction there have been results which characterise $M_n$ modulo primes $2$, $3$, and $5$. We collect these below and then discuss the similarity with results for larger primes.

\begin{thm}
\label{mod2}
(Theorem 5.5 of \cite{Eu2008}). The $n$th Motzkin number $M_n$ is even if and only if 
$$
\mbox{$n = (4i + \epsilon)4^{j+1} - \delta$ for $i, j \in \mathbb{N}, \epsilon \in \{1, 3\}$ and $\delta \in \{1, 2\}$. }
$$ 
\end{thm}
\bigskip

\begin{thm}
\label{mod3}
(Corollary 4.10 of \cite{Sagan2006}). Let T (\, 01 )\, be the set of numbers which have a base-$3$ representation consisting of the digits $0$ and $1$ only. Then the Motzkin numbers satisfy
\begin{align*}
M_n &\equiv  -1 \mod 3 \quad if  \quad n \in 3T (\, 01 )\, - 1, \\
M_n &\equiv 1 \mod 3  \quad if  \quad n \in 3T(\, 01 ) \quad or \quad n \in 3T(\, 01 )\, - 2, \\
M_n &\equiv 0 \mod 3 \quad otherwise.
\end{align*}
\end{thm}

\bigskip

\begin{thm}
\label{mod5}
(Theorem 5.4 of \cite{Sagan2006}). The Motzkin number $M_n$ is divisible by $5$ if and only if $n$ is one of the following forms
$$
(\, 5i + 1 )\, 5^{2j}  - 2,\, (\, 5i + 2 )\, 5^{2j-1}  - 1,\, (\, 5i + 3 )\, 5^{2j-1}  - 2,\,  (\, 5i + 4 )\, 5^{2j}  - 1
$$
where $i, j \in \mathbb{N}$ and $j \geq 1$.
\end{thm}

\bigskip
The above results and others have been used to establish asymptotic densities of the sets
\begin{equation}
\mbox{$S_q(0) = \{ n < N: M_n \equiv 0 \mod q \}$ .}
\end{equation}
for $q = 2, 4, 8, 3$ and $5$ - see \cite{RY2013} and \cite{Burns:2016vo}. In particular the asymptotic density of $S_2(0)$ is $\frac{1}{3}$ (\cite{RY2013} example 3.12), the asymptotic density of $S_4(0)$ is $\frac{1}{6}$ (\cite{RY2013} example 3.14), the asymptotic density of $S_3(0)$ is $1$ \cite{Burns:2016vo} and the asymptotic density of $S_5(0)$ is $\frac{1}{10}$ \cite{Burns:2016vo}.

There are 2 questions which we will investigate in this article. Firstly, what is the asymptotic density of $S_p(0)$ for primes from $7$ to $29$? Secondly, what structural features are evident in the distribution of $M_n$ modulo these primes? The investigation will proceed by constructing an automaton for each prime following the instructions from \cite{RY2013} (Algorithm 1). The state diagram for the automaton provides an excellent tool for analysing the behaviour of $M_n$ modulo each prime. It will be seen that the primes considered here either fit the pattern of theorems~\ref{mod2} or \ref{mod5} in which case the asymptotic density of $S_p(0)$ is $\frac{2}{p(p-1)}$ (or thereabouts), or fit the pattern of theorem~\ref{mod3} in which case the asymptotic density of $S_p(0)$ is $1$. For those primes which fit the pattern of theorems~\ref{mod2} and \ref{mod5} we will be interested in sets of the form
\begin{equation}
\label{S}
S(\, q, r, s, t )\, = \{ (\, qi + r )\, q^{sj + t} : i, j \in \mathbb{N} \}
\end{equation}
and
\begin{equation}
\label{S'}
S^{'}(\, q, r, s, t )\, = \{ (\, qi + r )\, q^{sj + t} : i, j \in \mathbb{N}, j \geq 1 \}
\end{equation}
for integers $q, r, s, t$.
\bigskip
For these sets we have the following theorem from \cite{Burns:2016vo}
\begin{thm}
\label{asyden}
Let $q, r, s, t \in \mathbb{Z} $ with $q, s > 0$, $t \geq 0$ and $0 \leq r < q$. Then the asymptotic density of the set $S$ is $(\, q^{t + 1 - s} (\, q^{s} - 1 )\, )\, ^{-1} $.  The asymptotic density of the set $S^{'}$ is \mbox{$(\, q^{t + 1} (\, q^{s} - 1 )\, )\, ^{-1} $}.
\end{thm}

\bigskip

\section{Background on automata}
Rowland and Yassawi showed in \cite{RY2013} that the behaviour of sequences such as $M_n~\mod~p$ can be studied by the use of finite state automata. The automaton has a finite number of states and rules for transitioning from one state to another. In the form described in \cite{RY2013} each state $s$ is represented by a polynomial in 2 variables $x$ and $y$. Each state has a value obtained by evaluating the polynomial at $x = 0$ and $y = 0$. In summary each $n$ is first represented in base $p$. The base $p$ digits of $n$ are fed into the automaton starting with the least significant digit. The automaton starts at an initial state and transitions to a new state as each digit is fed into it. The value of the final state after all $n$'s digits have been used is equal to $M_n~mod~p$. Refer to \cite{RY2013} for more details. 

In the remainder of this article we will provide details of the automaton for each prime. We will provide the polynomials and values for the states but not all of the transitions as these can be quite numerous. States are listed as $s[1]$, $s[2]$, ...Transitions, when provided, will be in the form $(\, i, j)\,  \to k$ which means that if the automaton is in state $i$ and receives digit $j$ then it will move to state $k$.  We will call a state $s$ a {\em \bf loop} state if all transitions from $s$ go to $s$ itself, i.e. \mbox{$(\, s, j)\,  \to s$} for all choices of $j$.

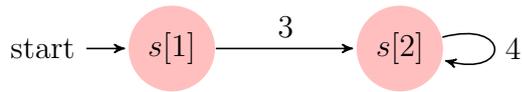
\begin{figure}[tbp]
\begin{tikzpicture}[->,>=stealth',shorten >=1pt,auto,node distance=3.0cm,
                    semithick]
  \tikzstyle{every state}=[fill=pink,draw=none,text=black]

  \node[initial,state] (A)                    {$s[1]$};
  \node[state]         (B) [right of=A] {$s[2]$};

  \path (A) edge    node {3} (B)
        (B) edge [loop right] node {4} (B);                      
 \end{tikzpicture}
                       \caption{Example of transition from state s[1] to state s[2] and a loop.}
                        \label{transition}
\end{figure}

States and transitions are represented visually in the form of a directed graph. For example, figure \ref{transition} represents an automaton which moves from state $s[1]$ to state $s[2]$ when it receives the digit $3$. It also moves from state $s[2]$ to state $s[2]$ (i.e. loops) if it is in state $s[2]$ and receives a digit $4$.

\bigskip

\section {Motzkin numbers modulo 7}
The automaton for $M_n~\mod 7$ has 11 states. These are
\begin{align*}
s[1] &= (5*y^4 + 5*y^3)*x^2 + 6*y^2*x + y \\
s[2] &= (2*y^3 + 2*y^2)*x^2 + y*x \\
s[3] &= 1 \\
s[4] &= 2 \\
s[5] &= 4 \\
s[6] &= (6*y^2 + 6*y)*x \\
s[7] &= (y^2 + y)*x + 2 \\
s[8] &= 3 \\
s[9] &= 0 \\
s[10] &= 5 \\
s[11] &= 6 \\
\end{align*}
The value of each state is obtained by setting $(\, x, y \,)\, = (\, 0, 0 \,)\,$.  The states with value $0$ are $s[1]$, $s[2]$, $s[6]$ and $s[9]$. States $s[1]$ and $s[2]$ cannot be final states because no transition ends at $s[1]$ and the only way to get to $s[2]$ is from $s[1]$ with a $0$ digit and $0$ never appears as a final (i.e. most significant) digit of a number. State $s[9]$ is the only loop state. The states of the automaton can be partitioned into $2$ sets. Let the set $A$ be defined by
$$
A := \{ \, s[1], s[6], s[7] \, \}
$$
and the set $\overline{A}$ be the set of all other states. If the state path is in a state from $\overline{A}$ and receives a $3$ then it transitions to state $s[9]$ from where it never leaves. Once the state path enters a state from the set $\overline{A}$ it never returns to any element of the set $A$. If the state path is in a state from $A$ and receives a $3$ then it transitions to a state from $\overline{A}$. Putting all this together, any number $n$ which has a base-$7$ representation containing $2$ or more digits $3$ will have a state path that ends at state $s[9]$, implying that \mbox{$M_n \equiv 0 \mod 7$}. Since the asymptotic density of the set of numbers having $0$ or $1$ digits $3$ in its base-$7$ representation is $0$, this implies that \textbf{the set} $\mathbf{S_{7}(0)}$ \textbf{has asymptotic density} $\mathbf{1}.$

\bigskip

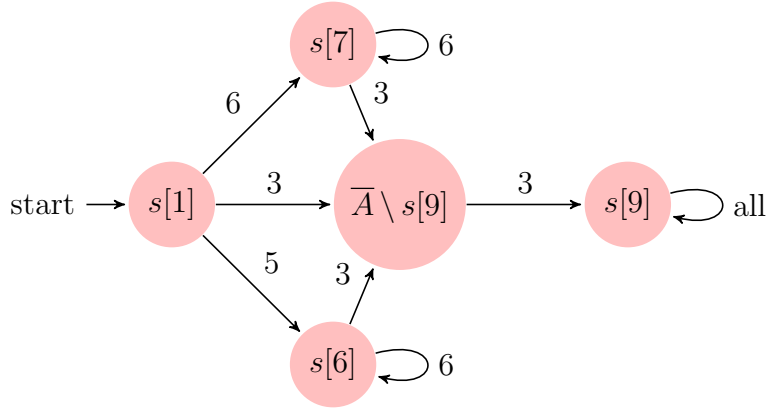
\begin{figure}[tbp]
\begin{tikzpicture}[->,>=stealth',shorten >=1pt,auto,node distance=3.0cm,
                    semithick]
  \tikzstyle{every state}=[fill=pink,draw=none,text=black]

  \node[initial,state] (A)                    {$s[1]$};
  \node[state]         (B) [below right of=A] {$s[6]$};
  \node[state]         (C) [right of=A] {${\overline{A}} \, \backslash \, s[9]$};
  \node[state]         (D) [above right of=A] {$s[7]$};
  \node[state]         (E) [right of=C]       {$s[9]$};

  \path (A) edge    node {5} (B)
            edge node {3} (C)
            edge node {6} (D)
        (B) edge [loop right] node {6} (B)
        edge node {3} (C)
         (C) edge node {3} (E)
                      (D) edge [loop right] node {6} (D)
                      edge node {3} (C)
                      (E) edge [loop right] node {all} (E); 
                      
 \end{tikzpicture}
                       \caption{Partial state diagram for $M_n \mod 7$.}
                        \label{mod7figure}
\end{figure}

Numbers $n$ of the form $7^k - 2$ for some $k~\geq~1$ have base-$7$ representation
$$
[\, n ]\,_{7} = \langle 6, 6, 6, ... 6, 5 \rangle
$$ 
where there are $(k-1)\,$ $6$'s in the representation. The subset of the state diagram given in figure~\ref{mod7figure} shows that the state path for this $n$ goes from the initial state to state $s[6]$ and remains there. Since the vale of $s[6]$ is $0$ this shows that
\begin{equation}
\label{7k-2mod7}
\mbox{$M_{n} \equiv 0 \mod 7 \,$ when $n = 7^k - 2$.}
\end{equation}

Similarly, when $n=7^k-1$ for some $k~\geq~1$, n has base-$7$ representation
$$
[\, n ]\,_{7} = \langle 6, 6, 6, ... 6, 6 \rangle
$$ 
where there are $k$ $6$'s in the representation. In this case the state path goes from the initial state to state $s[7]$ and stays there. Since the value of $s[7]$ is $2$ we have
\begin{equation}
\label{7k-1mod7}
\mbox{$M_{n} \equiv 2 \mod 7 \,$ when $n = 7^k - 1$.}
\end{equation}

\bigskip

Figure~\ref{mod7figure} provides a pictorial summary of the situation for $M_n~\mod~7$.

\bigskip

\section {Motzkin numbers modulo 11}
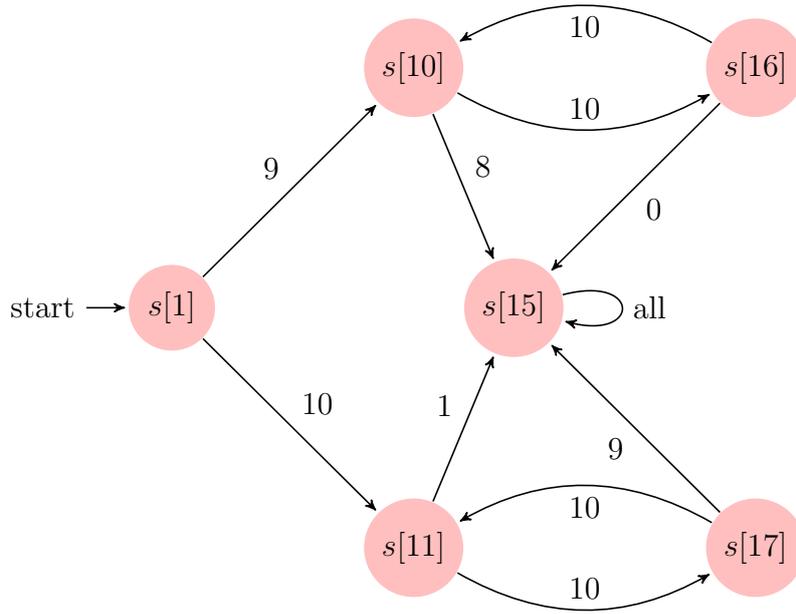
\begin{figure}[tbp]
\begin{tikzpicture}[->,>=stealth',shorten >=1pt,auto,node distance=4.5cm,
                    semithick]
  \tikzstyle{every state}=[fill=pink,draw=none,text=black]

  \node[initial,state] (A)                    {$s[1]$};
  \node[state]         (B) [below right of=A] {$s[11]$};
    \node[state]         (C) [right of=A] {$s[15]$};
  \node[state]         (D) [above right of=A] {$s[10]$};
  \node[state]         (E) [right of=D]       {$s[16]$};
  \node[state]         (F) [right of=B]       {$s[17]$};

  \path (A) edge     node {10} (B)
            edge node {9} (D)
        (B) edge [bend right] node {10} (F)
        edge node {1} (C)
         (C) edge [loop right] node {all} (C)
           (D) edge [bend right] node {10} (E)
                edge node {8} (C)
                      (E) edge [bend right] node {10} (D)
                      edge node {0} (C)
                       (F) edge [bend right] node {10} (B)
                      edge node {9} (C);
 \end{tikzpicture}
                       \caption{Partial state diagram for $M_n \mod 11$.}
                        \label{mod11figure}
\end{figure}

The automaton for $M_n~\mod 11$ has 17 states. These are
\begin{align*}
s[1] &= (9*y^4 + 9*y^3)*x^2 + 10*y^2*x + y \\
s[2] &= (2*y^3 + 2*y^2)*x^2 + y*x \\
s[3] &= 1 \\
s[4] &= 2 \\
s[5] &= 4 \\
s[6] &= 9 \\
s[7] &= 10 \\
s[8] &= 7 \\
s[9] &= 6 \\
s[10] &= (10*y^2 + 10*y)*x + 10 \\
s[11] &= (y^2 + y)*x + 10 \\
s[12] &= 3 \\
s[13] &= 8 \\
s[14] &= 5 \\
s[15] &= 0 \\
s[16] &= (10*y^2 + 10*y)*x \\
s[17] &= (y^2 + y)*x + 2 \\
\end{align*}
The states with value $0$ are $s[1]$, $s[2]$, $s[15]$ and $s[16]$. Again states $s[1]$ and $s[2]$ cannot be final states because the only way to get to $s[2]$ is from $s[1]$ with a $0$ digit and $0$ never appears as a final (i.e. most significant) digit of a number. State $s[15]$ is the only loop state. Figure \ref{mod11figure} is a subgraph of the state graph of this automaton which shows all the ways in which a state path can reach states $s[15]$ or $s[16]$. Figure~\ref{mod11figure} demonstrates that there are only 4 forms of numbers which reach a state with value $0$. These forms can be written either in terms of their base-$11$ representation or via formulae similar to those in theorems~\ref{mod2} and \ref{mod5}. In terms of their base-$11$ representations the $4$ forms are
$$
\mbox{$[\, n  \,]\,_{11} = \langle \, [\, i \,]\,_{11}, 8, 10, 10, ..., 10, 9 \rangle \,$ with an even number of $10$'s}
$$
$$
\mbox{$[\, n  \,]\,_{11} = \langle \, [\, i \,]\,_{11}, 0, 10, 10, ..., 10, 9 \rangle \,$ with an odd number of $10$'s}
$$
$$
\mbox{$[\, n \,]\,_{11} = \langle \, [\, i \,]\,_{11}, 1, 10, 10, ..., 10, 10 \rangle \,$ with an odd number of $10$'s}
$$
$$
\mbox{$[\, n  \,]\,_{11} = \langle [\, i \,]\,_{11}, 9, 10, 10, ..., 10, 10 \rangle \,$ with an even number of $10$'s}
$$
where $i$ is an arbitrary integer. Converting these forms to a formula we have
\bigskip
\begin{thm}
\label{mod11thm}
The Motzkin number $M_n$ is divisible by $11$ if and only if $n$ is one of the following forms
$$
(\, 11i + 9 )\, 11^{2j - 1}  - 2,\, (\, 11i + 1 )\, 11^{2j}  - 2,\, (\, 11i + 2 )\, 11^{2j-1}  - 1,\,  (\, 11i + 10 )\, 11^{2j}  - 1
$$
where $i, j \in \mathbb{N}$ and $j \geq 1$.
\end{thm}

\bigskip

The asymptotic densities of these 4 forms can be derived from theorem~\ref{asyden}.
\bigskip
\begin{thm}
\label{mod11ad}
The asymptotic density of numbers of the first and third forms in theorem~\ref{mod11thm} is $\frac {1}{120}$. The asymptotic density of numbers of the second and fourth forms in theorem~\ref{mod11thm} is $\frac {1}{1320}$ each. Therefore the asymptotic density of $S_{11}(0)$ is $\frac{1}{55}$.
\end{thm}

\bigskip

Figure~\ref{mod11figure} provides a pictorial summary of the situation for $M_n~\mod~11$.

\bigskip

\section {Motzkin numbers modulo 13}
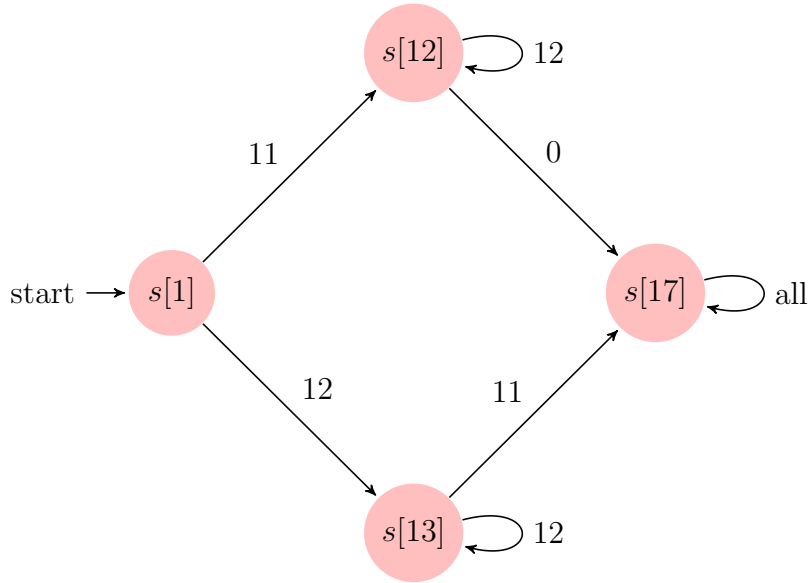
\begin{figure}[tbp]
\begin{tikzpicture}[->,>=stealth',shorten >=1pt,auto,node distance=4.5cm,
                    semithick]
  \tikzstyle{every state}=[fill=pink,draw=none,text=black]

  \node[initial,state] (A)                    {$s[1]$};
  \node[state]         (B) [above right of=A] {$s[12]$};
    \node[state]         (C) [below right of=A] {$s[13]$};
  \node[state]         (D) [above right of=C] {$s[17]$};

  \path (A) edge     node {11} (B)
            edge node {12} (C)
        (B) edge [loop right] node {12} (B)
        edge node {0} (D)
         (C) edge [loop right] node {12} (C)
               edge node {11} (D)
           (D) edge [loop right] node {all} (D);
 \end{tikzpicture}
                       \caption{Partial state diagram for $M_n \mod 13$.}
                        \label{mod13figure}
\end{figure}

The automaton for $M_n~\mod 13$ has 17 states. These are
\begin{align*}
s[1] &= (11*y^4 + 11*y^3)*x^2 + 12*y^2*x + y \\
s[2] &= (2*y^3 + 2*y^2)*x^2 + y*x \\
s[3] &= 1 \\
s[4] &= 2 \\
s[5] &= 4 \\
s[6] &= 9 \\
s[7] &= 8 \\
s[8] &= 12 \\
s[9] &= 10 \\
s[10] &= 11 \\
s[11] &= 3 \\
s[12] &= (12*y^2 + 12*y)*x \\
s[13] &= (y^2 + y)*x + 2 \\
s[14] &= 7 \\
s[15] &= 6 \\
s[16] &= 5 \\
s[17] &= 0. \\
\end{align*}
The states with value $0$ are $s[1]$, $s[2]$, $s[12]$ and $s[17]$. We exclude states $s[1]$ and $s[2]$ again. State $s[17]$ is the only loop state. Figure \ref{mod13figure} is a subgraph of the state graph of this automaton which shows all the ways in which a state path can reach states $s[12]$ or $s[17]$. Figure~\ref{mod13figure} demonstrates that there are only 2 forms of numbers which reach a state with value $0$. These forms can be written either in terms of their base-$13$ representation or via formulae similar to those in theorems~\ref{mod2} and \ref{mod5}. In terms of their base-$13$ representations the $2$ forms are
$$
\mbox{$[\, n  \,]\,_{13} = \langle \, [\, i \,]\,_{13}, 0, 12, 12, ..., 12, 11 \rangle \,$} 
$$
$$
\mbox{$[\, n  \,]\,_{13} = \langle \, [\, i \,]\,_{13}, 11, 12, 12, ..., 12, 12 \rangle \,$} 
$$
where $i$ is an arbitrary integer. Converting these forms to a formula we have
\bigskip

\begin{thm}
\label{mod13thm}
The Motzkin number $M_n$ is divisible by $13$ if and only if $n$ is one of the following forms
$$
(\, 13i + 1 )\, 13^{k}  - 2,\, (\, 13i + 12 )\, 13^{k}  - 1
$$
where $k \in \mathbb{N}$ and $k \geq 1$.
\end{thm}

\bigskip

The asymptotic densities of these 2 forms can be derived from theorem~\ref{asyden}.
\bigskip
\begin{thm}
\label{mod13ad}
The asymptotic density of each of the 2 sets of numbers described in theorem~\ref{mod13thm} is $\frac {1}{156}$. Therefore the asymptotic density of $S_{13}(0)$ is $\frac{1}{78}$.
\end{thm}

\bigskip

Figure~\ref{mod13figure} provides a pictorial summary of the situation for $M_n~\mod~13$.

\section {Motzkin numbers modulo 17}
The automaton for $M_n \mod 17$ has 23 states. These are
\begin{align*}
s[1] &= (15*y^4 + 15*y^3)*x^2 + 16*y^2*x + y \\
s[2] &= (2*y^3 + 2*y^2)*x^2 + y*x  \\
s[3] &= 1 \\
s[4] &= 2 \\
s[5] &= 4 \\
s[6] &= 9 \\
s[7] &= 0 \\
s[8] &= 8 \\
s[9] &= 12 \\
s[10] &= 7 \\
s[11] &= 15 \\
s[12] &= 6 \\
s[13] &= (16*y^2 + 16*y)*x + 16 \\
s[14] &= (y^2 + y)*x + 16 \\
s[15] &= 3 \\
s[16] &= 5 \\
s[17] &= 11 \\
s[18] &= 16 \\
s[19] &= 14 \\
s[20] &= 10 \\
s[21] &= 13 \\
s[22] &= (16*y^2 + 16*y)*x \\
s[23] &= (y^2 + y)*x + 2 \\
\end{align*}
The states with value $0$ are $s[1]$, $s[2]$, $s[7]$ and $s[22]$. We ignore states $s[1]$ and $s[2]$ again. State $s[7]$ is the only loop state. The states of the automaton can be partitioned into $2$ sets. Let the set $A$ be defined by
$$
A := \{ \, s[1], s[13], s[14], s[22], s[23] \, \}
$$
and the set $\overline{A}$ be the set of all other states. If the state path is in a state from $\overline{A}$ and receives a $5$ or $11$ then it transitions to state $s[7]$ from where it never leaves. Once the state path enters a state from the set $\overline{A}$ it never returns to any element of the set $A$. If the state path is in a state from $A$ and receives a $5$ or $11$ then it transitions to a state from $\overline{A}$. Putting all this together, any number $n$ which has a base-$17$ representation containing $2$ or more digits from $\{ 5, 11 \}$ will have a state path that ends at state $s[7]$, implying that \mbox{$M_n \equiv 0 \mod 17$}. Since the asymptotic density of the set of numbers having $0$ or $1$ digits from $\{ 5, 11 \}$ in its base-$17$ representation is $0$, this implies that \textbf{the set} $\mathbf{S_{17}(0)}$ \textbf{has asymptotic density} $\mathbf{1}.$

\bigskip
In the state diagram there is a cycle between states $s[13]$ and  $s[22]$ as follows:
\begin{align*}
(\,  s[13], 16 \, )\, &\to s[22] \\
(\,   s[22], 16 \, )\, &\to s[13] \\
\end{align*}
We also have the transition
$$
( \, s[1], 15 \, ) \to s[13].
$$
So numbers which have a base 17 representation of the form
$$
\mbox{$[\, n  \,]\,_{17} = \langle \, 16, 16, ..., 16, 15 \rangle \,$} 
$$
will have $M_n \equiv 0$ or $16 \mod 17$ depending on whether the number of $16$ digits is odd or even. In particular,
\begin{equation}
\mbox{if $\, n = 17^{2k} - 2 \,$ then $ \, M_n \equiv 0 \mod 17$.}
\end{equation}
and
\begin{equation}
\mbox{ if $\, n = 17^{2k-1} - 2 \,$ then $ \, M_n \equiv -1 \mod 17$.}
\end{equation}

\bigskip

In the state diagram there is also a cycle between states $s[14]$ and $s[23]$ and a similar argument shows that numbers which have a base $17$ representation of the form
$$
\mbox{$[\, n  \,]\,_{17} = \langle \, 16, 16, ..., 16, 16 \rangle \,$} 
$$
will have $M_n \equiv 2$ or $16 \mod 17$ depending on whether the number of $16$ digits is even or odd. In particular,
\begin{equation}
\mbox{if $n = 17^{2k} - 1$ then $M_n \equiv 2 \mod 17$} 
\end{equation}
and
\begin{equation}
\mbox{ if $\, n = 17^{2k-1} - 1 \,$ then $ \, M_n \equiv -1 \mod 17$.}
\end{equation}

\bigskip

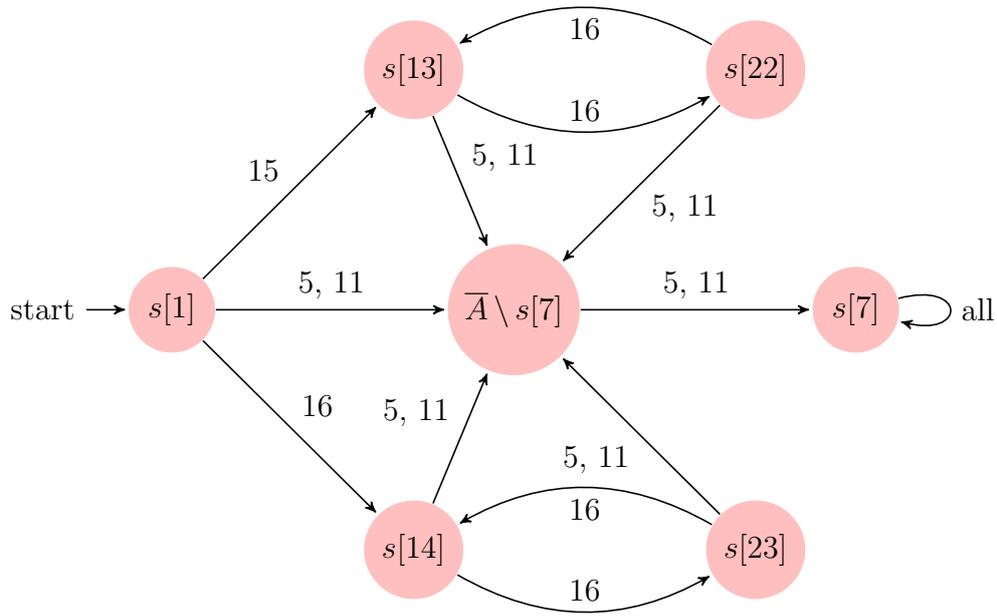
\begin{figure}[tbp]
\begin{tikzpicture}[->,>=stealth',shorten >=1pt,auto,node distance=4.5cm,
                    semithick]
  \tikzstyle{every state}=[fill=pink,draw=none,text=black]

  \node[initial,state] (A)                    {$s[1]$};
  \node[state]         (B) [below right of=A] {$s[14]$};
  \node[state]         (C) [right of=A] {${\overline{A}} \, \backslash \, s[7]$};
  \node[state]         (D) [above right of=A] {$s[13]$};
  \node[state]         (E) [right of=C]       {$s[7]$};
  \node[state]         (F) [right of=D]       {$s[22]$};
  \node[state]         (G) [right of=B]       {$s[23]$};

  \path (A) edge              node {16} (B)
            edge node {5, 11} (C)
            edge node {15} (D)
        (B) edge [bend right] node {16} (G)
        edge node {5, 11} (C)
         (C) edge node {5, 11} (E)
           (D) edge [bend right] node {16} (F)
                edge node {5, 11} (C)
                      (E) edge [loop right] node {all} (E)
                      (F) edge [bend right] node {16} (D)
                      edge node {5, 11} (C)
                       (G) edge [bend right] node {16} (B)
                      edge node {5, 11} (C);
 \end{tikzpicture}
                       \caption{Partial state diagram for $M_n \mod 17$.}
                        \label{mod17figure}
\end{figure}

Figure~\ref{mod17figure} provides a pictorial summary of the situation for $M_n~\mod~17$.

\bigskip

\section {Motzkin numbers modulo 19}
The automaton for $M_n \mod 19$ has 23 states. These are
\begin{align*}
s[1] &= (17*y^4 + 17*y^3)*x^2 + 18*y^2*x + y \\
s[2] &= (2*y^3 + 2*y^2)*x^2 + y*x \\
s[3] &= 1 \\
s[4] &= 2 \\ 
s[5] &= 4 \\
s[6] &= 9 \\
s[7] &= 13 \\
s[8] &= 0 \\ 
s[9] &= 18 \\
s[10] &= 3 \\ 
s[11] &= 7 \\ 
s[12] &= 16 \\ 
s[13] &= 14 \\ 
s[14] &= 17 \\
s[15] &= 6 \\ 
s[16] &= (18*y^2 + 18*y)*x \\
s[17] &= (y^2 + y)*x + 2 \\
s[18] &= 8 \\
s[19] &= 5 \\
s[20] &= 12 \\
s[21] &= 11 \\
s[22] &= 10 \\
s[23] &= 15 \\
\end{align*}
The states with value $0$ are $s[1]$, $s[2]$, $s[8]$ and $s[16]$. We ignore states $s[1]$ and $s[2]$ again. State $s[8]$ is the only loop state. The states can be partitioned into $2$ sets. Let the set $A$ be defined by
$$
A := \{ \, s[1], s[16], s[17] \, \}
$$
and the set $\overline{A}$ be the set of all other states. If the state path is in a state from $\overline{A}$ and receives a $4$ or $14$ then it transitions to state $s[8]$ from where it never leaves. Once the state path enters a state from the set $\overline{A}$ it never returns to any element of the set $A$. If the state path is in a state from $A$ and receives a $4$ or $14$ then it transitions to a state from $\overline{A}$. Putting all this together, any number $n$ which has a base-$19$ representation containing $2$ or more digits from $\{ 4, 14 \}$ will have a state path that ends at state $s[8]$, implying that \mbox{$M_n \equiv 0 \mod 19$}. This implies that \textbf{the set} $\mathbf{S_{19}(0)}$ \textbf{has asymptotic density} $\mathbf{1}.$

\bigskip
The state diagram also shows that
\begin{equation}
\mbox{if $\, n = 19^k - 2 \,$ then $ \, M_n \equiv 0 \mod 19$.}
\end{equation}
and
\begin{equation}
\mbox{  if $ \, n = 19^k - 1 \,$ then $\, M_n \equiv 2 \mod 19$.}
\end{equation}

\bigskip

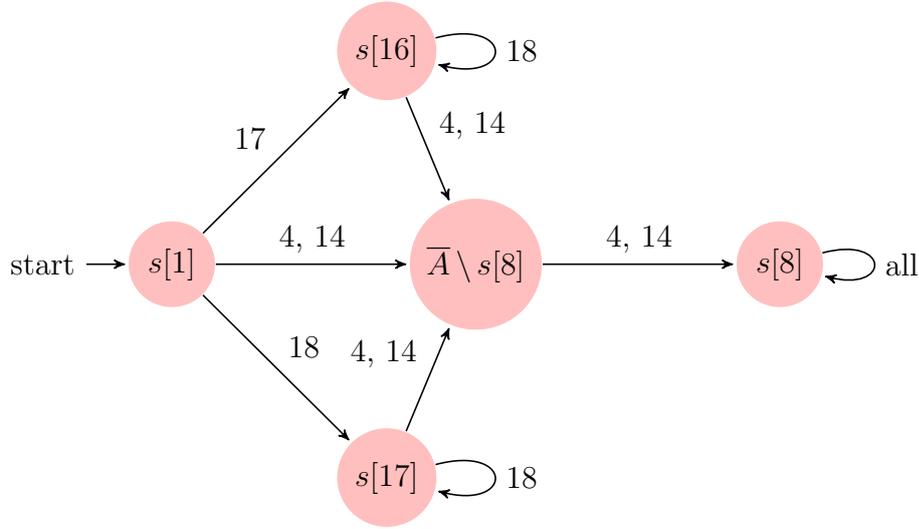
\begin{figure}[tbp]
\begin{tikzpicture}[->,>=stealth',shorten >=1pt,auto,node distance=4.0cm,
                    semithick]
  \tikzstyle{every state}=[fill=pink,draw=none,text=black]

  \node[initial,state] (A)                    {$s[1]$};
  \node[state]         (B) [below right of=A] {$s[17]$};
  \node[state]         (C) [right of=A] {${\overline{A}} \, \backslash \, s[8]$};
  \node[state]         (D) [above right of=A] {$s[16]$};
  \node[state]         (E) [right of=C]       {$s[8]$};

  \path (A) edge              node {18} (B)
            edge node {4, 14} (C)
            edge node {17} (D)
        (B) edge [loop right] node {18} (B)
        edge node {4, 14} (C)
         (C) edge node {4, 14} (E)
                      (D) edge [loop right] node {18} (D)
                      edge node {4, 14} (C)
                      (E) edge [loop right] node {all} (E); 
                      
 \end{tikzpicture}
                       \caption{Partial state diagram for $M_n \mod 19$.}
                        \label{mod19figure}
\end{figure}

Figure~\ref{mod19figure} provides a pictorial summary of the situation for $M_n~\mod~19$.

\bigskip

\section {Motzkin numbers modulo 23}
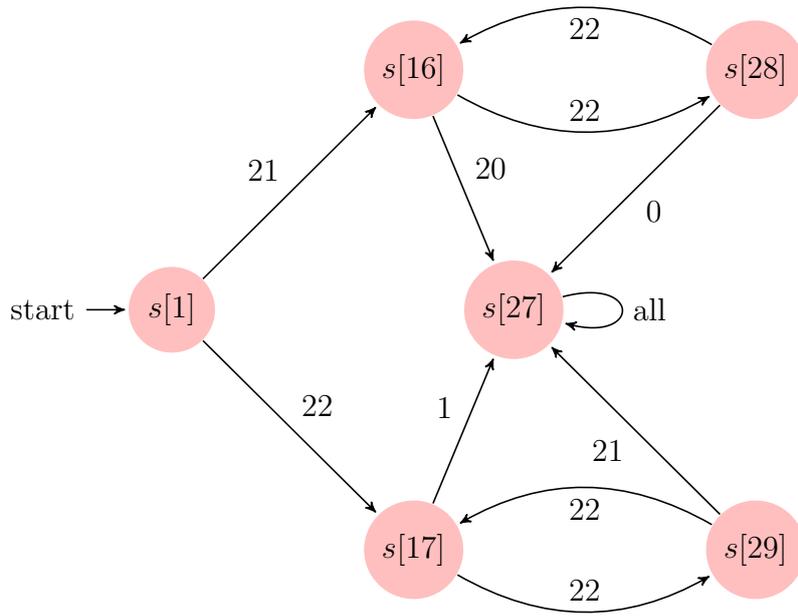
\begin{figure}[tbp]
\begin{tikzpicture}[->,>=stealth',shorten >=1pt,auto,node distance=4.5cm,
                    semithick]
  \tikzstyle{every state}=[fill=pink,draw=none,text=black]

  \node[initial,state] (A)                    {$s[1]$};
  \node[state]         (B) [below right of=A] {$s[17]$};
    \node[state]         (C) [right of=A] {$s[27]$};
  \node[state]         (D) [above right of=A] {$s[16]$};
  \node[state]         (E) [right of=D]       {$s[28]$};
  \node[state]         (F) [right of=B]       {$s[29]$};

  \path (A) edge     node {22} (B)
            edge node {21} (D)
        (B) edge [bend right] node {22} (F)
        edge node {1} (C)
         (C) edge [loop right] node {all} (C)
           (D) edge [bend right] node {22} (E)
                edge node {20} (C)
                      (E) edge [bend right] node {22} (D)
                      edge node {0} (C)
                       (F) edge [bend right] node {22} (B)
                      edge node {21} (C);
 \end{tikzpicture}
                       \caption{Partial state diagram for $M_n \mod 23$.}
                        \label{mod23figure}
\end{figure}

The automaton for $M_n~\mod 23$ has 29 states. These are
\begin{align*}
s[1] &= (21*y^4 + 21*y^3)*x^2 + 22*y^2*x + y \\
s[2] &= (2*y^3 + 2*y^2)*x^2 + y*x \\
s[3] &= 1 \\
s[4] &= 2 \\
s[5] &= 4 \\
s[6] &= 9 \\
s[7] &= 21 \\
s[8] &= 5 \\
s[9] &= 12 \\
s[10] &= 7 \\
s[11] &= 3 \\
s[12] &= 14 \\
s[13] &= 6 \\
s[14] &= 15 \\
s[15] &= 8 \\
s[16] &= (22*y^2 + 22*y)*x + 22 \\
s[17] &= (y^2 + y)*x + 22 \\
s[18] &= 19 \\
s[19] &= 11 \\
s[20] &= 22 \\
s[21] &= 10 \\
s[22] &= 16 \\
s[23] &= 20 \\
s[24] &= 17 \\
s[25] &= 18 \\
s[26] &= 13 \\
s[27] &= 0 \\
s[28] &= (22*y^2 + 22*y)*x \\
s[29] &= (y^2 + y)*x + 2 \\
\end{align*}
The states with value $0$ are $s[1]$, $s[2]$, $s[27]$ and $s[28]$. We exclude states $s[1]$ and $s[2]$ again. State $s[27]$ is the only loop state. Figure \ref{mod23figure} is a subgraph of the state graph of this automaton which shows all the ways in which a state path can reach states $s[27]$ or $s[28]$. Figure~\ref{mod23figure} demonstrates that there are only 4 forms of numbers which reach a state with value $0$. These forms can be written either in terms of their base-$23$ representation or via formulae similar to those in theorems~\ref{mod2} and \ref{mod5}. In terms of their base-$23$ representations the $4$ forms are
$$
\mbox{$[\, n  \,]\,_{23} = \langle \, [\, i \,]\,_{23}, 20, 22, 22, ..., 22, 21 \rangle \,$ with an even number of $22$'s}
$$
$$
\mbox{$[\, n  \,]\,_{23} = \langle \, [\, i \,]\,_{23}, 0, 22, 22, ..., 22, 21 \rangle \,$ with an odd number of $22$'s}
$$
$$
\mbox{$[\, n  \,]\,_{23} = \langle \, [\, i \,]\,_{23}, 1, 22, 22, ..., 22, 22 \rangle \,$ with an odd number of $22$'s}
$$
$$
\mbox{$[\, n  \,]\,_{23} = \langle \, [\, i \,]\,_{23}, 21, 22, 22, ..., 22, 22 \rangle \,$ with an even number of $22$'s}
$$
where $i$ is an arbitrary integer. Converting these forms to formulae we have
\bigskip
\begin{thm}
\label{mod23thm}
The Motzkin number $M_n$ is divisible by $23$ if and only if $n$ is one of the following forms
$$
(\, 23i + 21 )\, 23^{2j-1}  - 2,\, (\, 23i + 1 )\, 23^{2j}  - 2, (\, 23i + 2 )\, 23^{2j-1}  - 1,\, (\, 23i + 22 )\, 23^{2j}  - 1
$$
where $i, j \in \mathbb{N}$ and $j \geq 1$.
\end{thm}

\bigskip

The asymptotic densities of these 4 forms can be derived from theorem~\ref{asyden}.
\bigskip
\begin{thm}
\label{mod23ad}
The asymptotic density of numbers of the first and third forms in theorem~\ref{mod23thm} is $\frac {1}{528}$ each. The asymptotic density of numbers of the second and fourth forms in theorem~\ref{mod23thm} is $\frac {1}{12144}$ each. Therefore the asymptotic density of $S_{23}(0)$ is $\frac{1}{253}$.
\end{thm}

\bigskip

Figure~\ref{mod23figure} provides a pictorial summary of the situation for $M_n~\mod~23$.

\bigskip

\section {Motzkin numbers modulo 29}
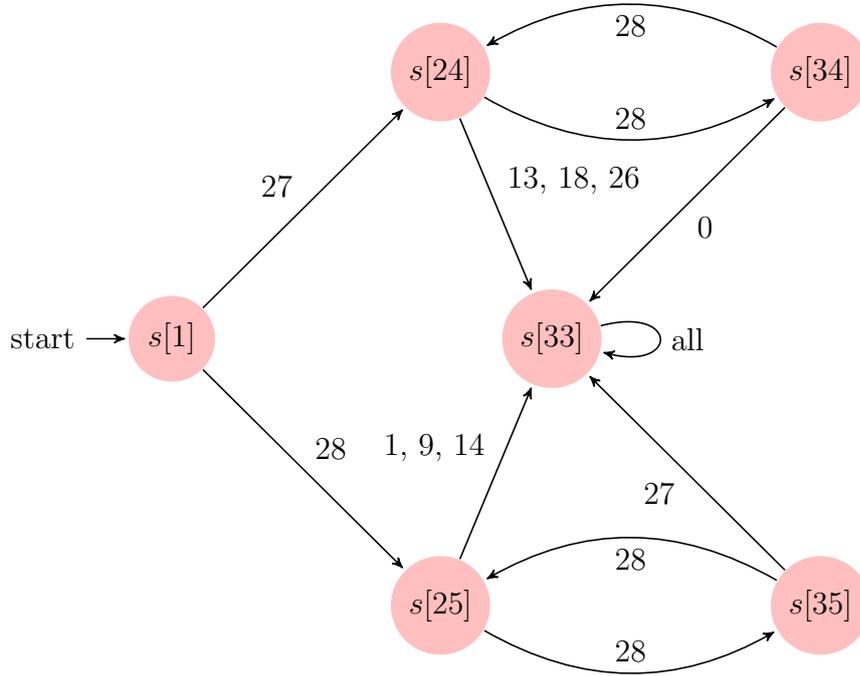
\begin{figure}[tbp]
\begin{tikzpicture}[->,>=stealth',shorten >=1pt,auto,node distance=5.0cm,
                    semithick]
  \tikzstyle{every state}=[fill=pink,draw=none,text=black]

  \node[initial,state] (A)                    {$s[1]$};
  \node[state]         (B) [below right of=A] {$s[25]$};
    \node[state]         (C) [right of=A] {$s[33]$};
  \node[state]         (D) [above right of=A] {$s[24]$};
  \node[state]         (E) [right of=D]       {$s[34]$};
  \node[state]         (F) [right of=B]       {$s[35]$};

  \path (A) edge     node {28} (B)
            edge node {27} (D)
        (B) edge [bend right] node {28} (F)
        edge node {1, 9, 14} (C)
         (C) edge [loop right] node {all} (C)
           (D) edge [bend right] node {28} (E)
                edge node {13, 18, 26} (C)
                      (E) edge [bend right] node {28} (D)
                      edge node {0} (C)
                       (F) edge [bend right] node {28} (B)
                      edge node {27} (C);
 \end{tikzpicture}
                       \caption{Partial state diagram for $M_n \mod 29$.}
                        \label{mod29figure}
\end{figure}

The automaton for $M_n~\mod 29$ has 35 states. These are
\begin{align*}
s[1] &= (27*y^4 + 27*y^3)*x^2 + 28*y^2*x + y \\
s[2] &= (2*y^3 + 2*y^2)*x^2 + y*x \\
s[3] &= 1 \\
s[4] &= 2 \\
s[5] &= 4 \\
s[6] &= 9 \\
s[7] &= 21 \\
s[8] &= 22 \\
s[9] &= 11 \\
s[10] &= 23 \\
s[11] &= 13 \\
s[12] &= 27 \\
s[13] &= 25 \\
s[14] &= 17 \\
s[15] &= 12 \\
s[16] &= 26 \\
s[17] &= 7 \\
s[18] &= 14 \\
s[19] &= 15 \\
s[20] &= 20 \\
s[21] &= 28 \\
s[22] &= 16 \\
s[23] &= 10 \\
s[24] &= (28*y^2 + 28*y)*x + 28 \\
s[25] &= (y^2 + y)*x + 28 \\
s[26] &= 3 \\
s[27] &= 19 \\
s[28] &= 5 \\
s[29] &= 6 \\
s[30] &= 18 \\
s[31] &= 8 \\
s[32] &= 24 \\
s[33] &= 0 \\
s[34] &= (28*y^2 + 28*y)*x \\
s[35] &= (y^2 + y)*x + 2 \\
\end{align*}
The states with value $0$ are $s[1]$, $s[2]$, $s[33]$ and $s[34]$. We exclude states $s[1]$ and $s[2]$ again. State $s[33]$ is the only loop state. Figure \ref{mod29figure} is a subgraph of the state graph of this automaton which shows all the ways in which a state path can reach states $s[33]$ or $s[34]$. Figure~\ref{mod29figure} demonstrates that there are only 4 forms of numbers which reach a state with value $0$. These forms can be written either in terms of their base-$29$ representation or via formulae similar to those in theorems~\ref{mod2} and \ref{mod5}. In terms of their base-$29$ representations the $4$ forms are
$$
\mbox{$[\, n  \,]\,_{29} = \langle \, [\, i \,]\,_{29}, a, 28, 28, ..., 28, 27 \rangle \,$ with an even number of $28$'s}
$$
$$
\mbox{$[\, n  \,]\,_{29} = \langle \, [\, i \,]\,_{29}, 0, 28, 28, ..., 28, 27 \rangle \,$ with an odd number of $28$'s}
$$
$$
\mbox{$[\, n  \,]\,_{29} = \langle \, [\, i \,]\,_{29}, b, 28, 28, ..., 28, 28 \rangle \,$ with an odd number of $28$'s}
$$
$$
\mbox{$[\, n  \,]\,_{29} = \langle \, [\, i \,]\,_{29}, 27, 28, 28, ..., 28, 28 \rangle \,$ with an even number of $28$'s}
$$
where $i$ is an arbitrary integer, $a \in \{ \, 13, 18, 26 \, \}$ and $b \in \{ \, 1, 9, 14 \, \}$. Converting these forms to formulae we have
\bigskip
\begin{thm}
\label{mod29thm}
The Motzkin number $M_n$ is divisible by $29$ if and only if $n$ is one of the following forms
$$
(\, 29i + a + 1 )\, 29^{2j-1}  - 2,\, (\, 29i + 1 )\, 29^{2j}  - 2, (\, 29i + b + 1 )\, 29^{2j-1}  - 1,\, (\, 29i + 28 )\, 29^{2j}  - 1
$$
where $a \in \{ \, 13, 18, 26 \, \}$, $b \in \{ \, 1, 9, 14 \, \}$$i, j \in \mathbb{N}$ and $j \geq 1$.
\end{thm}

\bigskip

The asymptotic densities of these forms can be derived from theorem~\ref{asyden}.
\bigskip
\begin{thm}
\label{mod29ad}
The asymptotic density of numbers of the first and third forms in theorem~\ref{mod29thm} is $\frac {1}{840}$ each. The asymptotic density of numbers of the second and fourth forms in theorem~\ref{mod29thm} is $\frac {1}{24360}$ each. Therefore the asymptotic density of $S_{29}(0)$ is $\frac{22}{3045}$.
\end{thm}

\bigskip

Figure~\ref{mod29figure} provides a pictorial summary of the situation for $M_n~\mod~29$.

\bigskip

\section{General observations}

\begin{observation}
The automaton for every prime we have looked at has a loop state and that state has value $0$. In addition, apart from the first $2$ states there were always $2$ other states with a value of $0$. We do not know of a theoretical reason why this should be the case but it would be interesting to either find a counterexample or prove that this always occurs. It should be noted that he number of states associated with the automaton depends on the algorithm used as mentioned in \cite{RY2013}.
\end{observation}

\bigskip

\begin{observation}
For primes $p$ up to $23$ the asymptotic density of $S_p(0)$ is either $1$ or $\frac{2}{p(p-1)}$. The result for $p~=~29$ shows that this is not true in general. However, it may still be the case that the asymptotic density of $S_p(0) \geq \frac{2}{p(p-1)}$ for all primes. If the state diagram for the automaton of $M_n~\mod~p$ contains a subgraph which looks like the one in figure \ref{p(p-1)figure}then the asymptotic density of $S_p(0) \geq \frac{2}{p(p-1)}$.  In figure \ref{p(p-1)figure} $w, x, y, z$ are arbitrary elements of $\frac{\mathbb{Z}}{p \mathbb{Z}}$. Figure~\ref{mod13figure} is in this form. Figures~\ref{mod11figure}, \ref{mod23figure} and \ref{mod29figure} are also in this form after an appropriate identification of nodes.

\begin{figure}[tbp]
\begin{tikzpicture}[->,>=stealth',shorten >=1pt,auto,node distance=4.5cm,
                    semithick]
  \tikzstyle{every state}=[fill=pink,draw=none,text=black]

  \node[initial,state] (A)                    {$A$};
  \node[state]         (B) [above right of=A] {$B$};
    \node[state]         (C) [below right of=A] {$C$};
  \node[state]         (D) [above right of=C] {$D$};

  \path (A) edge     node {w} (B)
            edge node {x} (C)
        (B) edge [loop right] node {p-1} (B)
        edge node {y} (D)
         (C) edge [loop right] node {p-1} (C)
               edge node {z} (D)
           (D) edge [loop right] node {all} (D);
 \end{tikzpicture}
                       \caption{Partial state diagram for $M_n \mod p$.}
                        \label{p(p-1)figure}
\end{figure}
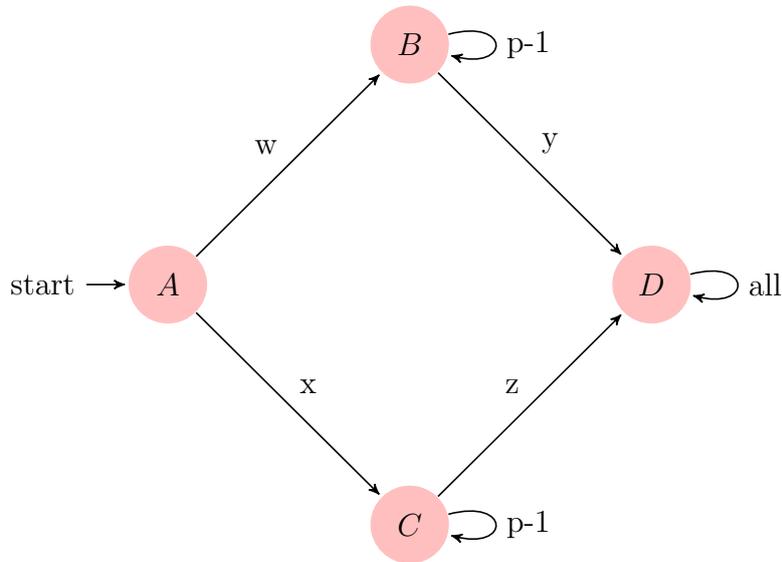

\end{observation}

\bigskip

\section{Acknowledgement}
We would like to thank Eric Rowland for introducing us to automata as a tool for examining Motzkin numbers.

\bigskip

\bibliographystyle{plain}
\begin{small}
\bibliography{ref}

\begin{thebibliography}{1}

\bibitem{Burns:2016vo}
Rob Burns.
\newblock Asymptotic density of motzkin numbers modulo small primes.
\newblock {\em ArXiv}, arXiv:1611.04910:6, 2016.

\bibitem{Sagan2006}
E.~Deutsch and B.E. Sagan.
\newblock Congruences for {C}atalan and {M}otzkin numbers and related
  sequences.
\newblock {\em Journal of Number Theory}, 117(1):191--215, 2006.

\bibitem{Eu2008}
Sen-Peng Eu, Shu-Chung Liu, and Yeong-Nan Yeh.
\newblock Catalan and {M}otzkin numbers modulo 4 and 8.
\newblock {\em European Journal of Combinatorics}, 29:1449--1466, 2008.

\bibitem{Krat2016}
C.~Krattenthaler and T.~W. M\"uller.
\newblock Motzkin numbers and related sequences modulo powers of $2$.
\newblock {\em ArXiv}, arXiv:1608.05657:28, 2016.

\bibitem{KM2013}
Christian Krattenthaler and Thomas~W. M\"uller.
\newblock A method for determining the mod-$3^k$ behaviour of recursive
  sequences.
\newblock {\em ArXiv}, arXiv:1308.2856:82, 2013.

\bibitem{RY2013}
Eric Rowland and Reem Yassawi.
\newblock Automatic congruences for diagonals of rational functions.
\newblock {\em ArXiv}, arXiv:1310.8635:42, 2013.

\end{thebibliography}
\end{small}

\end{document}